\documentclass[12pt,a4paper]{article}
\usepackage{amsmath}
\usepackage{amssymb}
\usepackage{amsthm}
\usepackage{amsfonts}
\usepackage{latexsym}

\theoremstyle{plain}
\newtheorem{thm}{Theorem}[section]
\newtheorem{cor}{Corollary}[section]
\newtheorem{lem}{Lemma}[section]

\theoremstyle{definition}

\theoremstyle{remark}

\title{A local Szeg\"{o}-type theorem\\ in Toeplitz quantization}
\author{Roberto Paoletti\footnote{\noindent{\bf Address:}
Dipartimento di Matematica e Applicazioni, Universit\`a degli Studi
di Milano Bicocca, Via R. Cozzi 53, 20125 Milano,
Italy; {\bf e-mail}: roberto.paoletti@unimib.it }}
\date{}

\begin{document}

\maketitle

\begin{abstract}
A Szeg\"{o}-type theorem for Toeplitz operators was proved by Boutet de Monvel and
Guillemin for general Toeplitz structures. We give a local version of this result in the setting of positive line bundles
on compact symplectic manifolds.
\end{abstract}

\section{Introduction}

The goal of this paper is to establish a local version of 
a Szeg\"{o}-type theorem for
Toeplitz operators proved by Boutet de Monvel and Guillemin  in \cite{bg} (Theorem 6 in \S 1,
see also \S 13), which is the analogue of results of Weinstein
\cite{w} and Widom \cite{wi}
for pseudodifferential operators. While \cite{bg} deals with general Toeplitz structures, here we shall focus on the
special case of positive line bundles, which is natural in algebraic geometry and geometric quantization.

We shall adopt the point of view of \cite{p-weyl}, \cite{p-equiv-asymp}, \cite{p-local-trace},
where certain spectral asymptotics for Toeplitz operators
are given a local interpretation building on the approach to generalized Szeg\"{o} kernels in \cite{z},
\cite{bsz} and \cite{sz}, based on the microlocal
theory of \cite{bs}. Thus the natural context for the present discussion is the category of quantized compact symplectic
manifolds and the spaces of \lq almost holomorphic sections\rq \, described in \cite{bg} and \cite{sz}; for
ease of exposition, we shall confine ourselves to the complex projective case.

Our setting and normalization conventions are as follows. Let $M$ be connected complex d-dimensional projective manifold, and
$A$ an ample line bundle on it. Let $h$ be an Hermitian metric on $A$ such that the unique compatible connection
has curvature $\Theta=-2i\omega$, where $\omega$ is a K\"{a}hler form on $M$. Let $A^*$ be the dual line bundle and
$X\subseteq A^*$ the unit circle bundle, with projection $\pi:X\rightarrow M$.
Then the connection 1-form $\alpha$ on $X$ is a contact form, and we shall adopt
$d\mu_M=:(1/\mathrm{d}!)\,\omega^{\wedge \mathrm{d}}$ and $d\mu_X=:(1/2\pi)\alpha\wedge \pi^*(d\mu_M)$ as volume forms
on $M$ and $X$, respectively. Given these choices, we shall identify (generalized) functions, densities and
half-densities on $X$.

Let $H(X)\subseteq L^2(X)$ be the Hardy space, and let $H(X)=\bigoplus _kH_k(X)$ be the
decomposition of $H(X)$ into $S^1$-equivariant summands; thus $H_k(X)$ is naturally unitarily
isomorphic to the space $H^0\left(M,A^{\otimes k}\right)$
of global holomorphic sections of $A^{\otimes k}$. The orthogonal projector
$\Pi:L^2(X)\rightarrow H(X)$ extends to a continuous operator $\Pi:\mathcal{D}'(X)\rightarrow \mathcal{D}'(X)$,
and the distributional kernel of the latter is the Szeg\"{o} kernel studied in \cite{bs}.
Thus, $\Pi=\bigoplus _k\Pi_k$, where $\Pi_k:L^2(X)\rightarrow H_k(X)$ is the level-$k$
Szeg\"{o} kernel.

A Toeplitz operator on $X$ of order $\nu\in \mathbb{Z}$ is a composition $T=:\Pi\circ Q\circ \Pi$,
where $Q$ is a pseudodifferential operator of classical type and order $\nu$ on $X$ \cite{bg}; we can naturally view
$T$ as linear operator on $H(X)$. If in particular $T$ is $S^1$-invariant,
then it leaves every summand $H_k(X)$ invariant, inducing \lq equivariant\rq \,endomorphisms
$T_k:H_k(X)\rightarrow H_k(X)$. Our focus will be on $S^1$-invariant zeroth order self-adjoint Toeplitz operators
and the (local) spectral asymptotics as $k\rightarrow +\infty$ of their equivariant pieces.
A notable example is $T_f=\Pi\circ M_f\circ \Pi$, where $M_f$ is multiplication
by a real valued $f\in \mathcal{C}^\infty(M)$; under the unitary isomorphism $H_k(X)\cong H^0\left(M,A^{\otimes k}\right)$, $T_k$ corresponds
to the Berezin-Toeplitz quantization of the classical observable $f$.

Thus, let $T$ be an $S^1$-invariant zeroth order self-adjoint Toeplitz operator on $X$, and for
$k=1,2,\ldots$ let $T_k:H_k(X)\rightarrow H_k(X)$ be the self-adjoint endomorphism induced by restriction.
Let $\lambda_{kj}$, $j=1,\ldots,N_k=\dim\big(H_k(X)\big)$ be the eigenvalues of $T_k$ repeated according to
multiplicity, and suppose that $e_{kj}$ is an eigenfunction of $T_k$ for the eigenvalue $\lambda_{kj}$,
so chosen that $(e_{kj})_j$ is an orthonormal basis of $H_k(X)$.

We shall study the asymptotics of the following measures on $\mathbb{R}$:
$$
\mathcal{T}_{x,k}=:\sum_{j=1}^{N_k}\delta_{\lambda_{kj}}\,\big|e_{kj}(x)\big|^2\,\,\,\,\,\,\,\,(x\in X,\,\,k=1,2,\ldots);
$$
$\mathcal{T}_{x,k}$ does not depend on the choice of the $e_{kj}$'s. Furthermore, $\big|e_{kj}(x)\big|^2$ only depends
on $m=\pi(x)$, hence we may write $\mathcal{T}_{m,k}$ for $\mathcal{T}_{x,k}$.

Before stating our main result, let us recall (after \cite{bg})
that the connection 1-form generates a closed symplectic cone
$$
\Sigma=:\big\{(x,r\alpha_x):\,x\in X,r>0\big\}\subseteq T^*X\setminus\{0\},
$$
and that the (principal) \textsl{symbol} $\sigma_T:\Sigma\rightarrow \mathbb{C}$
of the Toeplitz operator $T=\Pi\circ Q\circ \Pi$ is simply the restriction to $\Sigma$ of the principal
symbol $\sigma_Q$ of the pseudodifferential operator $Q$; if $T$ is self-adjoint, then $\sigma_T$ is real valued.

The \textit{reduced symbol} of $T$ is then the $\mathcal{C}^\infty$ function $\varsigma_T(x)=:\sigma_T(x,\alpha_x)$
 on $X$; if $T$ is $S^1$-invariant then so is $\varsigma_T$, which may thus be viewed as a
$\mathcal{C}^\infty$ function on $M$.
\begin{thm}
\label{thm:maim}
Let $T$ be an $S^1$-invariant zeroth order self-adjoint Toeplitz operator on $X$ and let
$\chi\in \mathcal{S}(\mathbb{R})$ be function of rapid decrease. Then uniformly in $m\in M$ the following asymptotic expansion holds as $k\rightarrow +\infty$:
\begin{eqnarray*}
\langle\mathcal{T}_{m,k},\chi\rangle&\sim& \left(\frac k\pi\right)^\mathrm{d}\,
\left(\chi\big(\varsigma_T(m)\big)\,+\sum_{j\ge 1}k^{-j}P_j(\chi)\big(\varsigma_T(m)\big)\right),
\end{eqnarray*}
where $P_j$ is a differential operator of degree $2j$ depending on $m$.
\end{thm}

We obtain a global asymptotic expansion by integrating over $M$. More precisely,
for $k=1,2,\ldots$, define measures on $\mathbb{R}$ given by
$\mathcal{T}_{k}=:\sum_{j=1}^{N_k}\delta_{\lambda_{kj}}$.

\begin{cor}
For any $\chi\in \mathcal{S}(\mathbb{R})$, the following asymptotic expansion holds for $k\rightarrow +\infty$:
\begin{eqnarray*}
\langle\mathcal{T}_{k},\chi\rangle&\sim& \left(\frac k\pi\right)^\mathrm{d}\,\
\left(\int_M \chi\big(\varsigma_T(m)\big)d\mu_M(m)+\sum_{j\ge 1}k^{-j}C_j(\chi)\right).
\end{eqnarray*}
\end{cor}

In the present setting of positive line bundles, the Szeg\"{o}-type theorem of \cite{bg} is the following:

\begin{cor}
For any $\chi\in \mathcal{S}(\mathbb{R})$,
\begin{eqnarray*}
\lim _{k\rightarrow +\infty}\left(\frac \pi k\right)^{\mathrm{d}}
\langle\mathcal{T}_{k},\chi\rangle=
\int_M \chi\big(\varsigma_T(m)\big)d\mu_M(m).
\end{eqnarray*}
\end{cor}

\section{Proof of Theorem \ref{thm:maim}}

To begin with, we may reduce the proof to the case where $T$ is elliptic with everywhere
positive reduced symbol; recall that a Toeplitz operator is called \textit{elliptic} if
its symbol is nowhere vanishing, hence if $T$ is self-adjoint then it is elliptic if and only
if  $\varsigma_T$ has constant sign.

In fact, choose $C>0$ such that $\varsigma_T+C>0$, and define
$\widetilde{T}=:T+c\Pi$; then $\widetilde{T}$ is an elliptic zeroth order $S^1$-invariant
self-adjoint Toeplitz operator, with reduced symbol $\varsigma_{\widetilde{T}}=\varsigma_T+c>0$
and eigenvalues $\widetilde{\lambda}_{kj}=\lambda_{kj}+c$ with the same eigenfunctions $e_{kj}$.
Now suppose that the statement
has been proved for $\widetilde{T}$, and define $\chi_c(\tau)=:\chi(\tau-c)$. Then if $\widetilde{\mathcal{T}}_{m,k}$
is the analogue of $\mathcal{T}_{m,k}$ for $\widetilde{T}$ we have
\begin{eqnarray*}
\lefteqn{\langle \mathcal{T}_{m,k},\chi\rangle =\sum_{j=1}^{N_k}\big|e_{kj}(x)\big|^2\langle\delta_{\lambda_{kj}},\chi\rangle
=\sum_{j=1}^{N_k}\big|e_{kj}(x)\big|^2\langle\delta_{\lambda_{kj}+c},\chi_c\rangle }\\
&=&
\left<\widetilde{\mathcal{T}}_{m,k},\chi_c\right>\sim
\left(\frac k\pi\right)^\mathrm{d}\,\left(\chi_c\big(\widetilde{\varsigma}_T(m)\big)+\sum_{j\ge 1}k^{-j}C_j\right)\\
&=&\left(\frac k\pi\right)^\mathrm{d}\,\left(\chi\big(\varsigma_T(m)\big)+\sum_{j\ge 1}k^{-j}C_j\right).
\end{eqnarray*}

Let us then assume in the following that $\varsigma_T>0$.

We shall denote by either $\mathcal{F}(\gamma)$ or
$\widehat{\gamma}$ the Fourier transform of
$\gamma\in \mathcal{S}'(\mathbb{R})$ (thus
$\widehat{f}(\xi)=:(2\pi)^{-1/2}\int_{-\infty}^{+\infty} f(x)\,e^{-i\xi x}dx$ if $f\in L^1$);
for any $\lambda\in \mathbb{R}$, we have
$\delta_\lambda=\left(1/\sqrt{2\pi}\right)\,\mathcal{F}\left(e^{i\lambda (\cdot)}\right)$. Define $\chi_k(\tau)=:\chi(\tau/k)$, so that
$\langle\delta_\lambda,\chi\rangle=:\langle \delta_{k\lambda},\chi_k\rangle$,
$\forall\,\lambda\in \mathbb{R}$.
We have if $m=\pi(x)$:
\begin{eqnarray}
\label{eqn:fourier-transform}
\lefteqn{\langle \mathcal{T}_{m,k},\chi\rangle =\sum_{j=1}^{N_k}\big|e_{kj}(x)\big|^2\langle\delta_{\lambda_{kj}},\chi\rangle
=\sum_{j=1}^{N_k}\big|e_{kj}(x)\big|^2\langle\delta_{k\lambda_{kj}},\chi_k\rangle}\nonumber\\
&=&\frac{1}{\sqrt{2\pi}}\,\sum_{j=1}^{N_k}\big|e_{kj}(x)\big|^2
\left<\mathcal{F}\left(e^{ik\lambda_{kj} (\cdot)}\right),\chi_k\right>=\frac{1}{\sqrt{2\pi}}\cdot \sum_{j=1}^{N_k}\big|e_{kj}(x)\big|^2\left<e^{ik\lambda_{kj} (\cdot)},\widehat{\chi}_k\right>\nonumber\\
&=&\frac{1}{\sqrt{2\pi}}\cdot \int_{-\infty}^\infty\left(\sum_{j=1}^{N_k}\big|e_{kj}(x)\big|^2
e^{ik\lambda_{kj} \tau}\right)\,\widehat{\chi}_k(\tau)\,d\tau.
\end{eqnarray}

To interpret the expression within brackets, let us introduce the first order Toeplitz
operator $T'=:D\circ T$, where $D=:-i\,\partial/\partial \theta$ and $\partial/\partial \theta$
is the generator of the $S^1$-action. Then $T'$ is $S^1$-invariant and self-adjoint, and has eigenvalues
$\lambda'_{kj}=:k\,\lambda_{kj}$, with eigensections $e_{kj}$. It is furthermore elliptic, with everywhere positive
symbol $\sigma_{T'}(x,r\alpha_x)=r\,\varsigma_T(x)$.

By the theory in \S 2 and \S 12 of \cite{bg}, there exists an elliptic, self-adjoint $S^1$-invariant operator
$Q$ on $X$ satisfying the following conditions: $[\Pi,Q]=0$; $T'=\Pi\circ Q\circ \Pi$; the principal symbol $\sigma_Q$
of $Q$ is everywhere positive on $T^*X\setminus \{0\}$. In particular, $Q$ leaves $H(X)$ invariant, and $T'$ is the
restriction of $Q$ to $H(X)$. If $U(\tau)=:e^{i\tau Q}$ ($\tau\in \mathbb{R}$), then $U(\tau)\big(e_{kj}\big)=e^{ik\lambda_{kj}\tau}\cdot e_{kj}$.
The distributional kernel
of $U_k(\tau)=:U(\tau)\circ \Pi_k$ is then
$$
U_k(\tau)(x,y)=:\sum_{j=1}^{N_k}e^{ik\lambda_{kj}\tau}\,
e_{kj}(x)\cdot \overline{e_{kj}(y)}\,\,\,\,\,\,\,(x,y\in X);
$$
therefore, (\ref{eqn:fourier-transform}) may be rewritten
\begin{eqnarray}
\label{eqn:fourier-transform-1}
\langle \mathcal{T}_{m,k},\chi\rangle& =&
\frac{1}{\sqrt{2\pi}}\cdot \int_{-\infty}^\infty\widehat{\chi}_k(\tau)\,U_k(\tau)(x,x)\,d\tau\nonumber\\
&=&\frac{1}{2\pi}\cdot \int_{-\infty}^{+\infty}\int_{-\infty}^{+\infty}\chi(\lambda/ k)e^{-i\tau\lambda}
\,U_k(\tau)(x,x)\,d\lambda\,d\tau\nonumber\\
&=&\frac{k}{2\pi}\cdot \int_{-\infty}^{+\infty}\int_{-\infty}^{+\infty}\chi(\lambda)e^{-ik\tau\lambda}
\,U_k(\tau)(x,x)\,d\lambda\,d\tau,
\end{eqnarray}
where in the last equality the change of variable $\lambda\rightarrow k\lambda$ has been performed.

Next we choose a sufficiently small $\epsilon>0$, and
$\beta\in \mathcal{C}^\infty\big((-\epsilon,\epsilon)\big)$ such that $\beta\ge 0$ and
$\beta(\tau)=1$ for $-1/2<\tau<1/2$. Inserting the identity $1=\beta(\tau)+\big(1-\beta(\tau)\big)$ in (\ref{eqn:fourier-transform-1}),
we obtain
$$
\langle \mathcal{T}_{m,k},\chi\rangle=\langle \mathcal{T}_{m,k},\chi\rangle '+\langle \mathcal{T}_{m,k},\chi\rangle'',
$$ where in the former summand the integrand has been multiplied by $\beta$, and in the latter by $1-\beta$.

\begin{lem}
As $k\rightarrow +\infty$, we have
$\langle \mathcal{T}_{m,k},\chi\rangle''=O\left(k^{-\infty}\right)$.
\end{lem}

\textit{Proof.}
In view of (\ref{eqn:fourier-transform-1}), we have
\begin{eqnarray}
\label{eqn:fourier-transform-2}
 \langle \mathcal{T}_{m,k},\chi\rangle''=\frac{k}{\sqrt{2\pi}}\,
 \int_{-\infty}^{+\infty}\widehat{\chi}(k\tau)\,\big(1-\beta(\tau)\big)\,
\,U_k(\tau)(x,x)\,d\tau.
 \end{eqnarray}
 Now since $\widehat{\chi}$ is of rapid decrease and $1-\beta$ is bounded and supported where $|\tau|\ge \epsilon/2$,
 for every $N>0$ we have
 $$
 \left|\widehat{\chi}(k\tau)\,\big(1-\beta(\tau)\big)\right|\le \frac{C_N}{k^{2N}\left(\tau^2+\epsilon^2\right)^N}.
 $$
On the other hand, by the Tian-Zelditch asymptotic expansion one has $\big|U_k(\tau)(x,x)\big|\le C\,k^\mathrm{d}$ uniformly in $\tau\in \mathbb{R}$ and
$x\in X$.

\hfill Q.E.D.

\bigskip

We thus need only consider the asymptotics of $\langle \mathcal{T}_{m,k},\chi\rangle'$. To proceed further,
let us rewrite (\ref{eqn:fourier-transform-1}) as
\begin{eqnarray}
\label{eqn:fourier-transform-comp}
\lefteqn{\langle \mathcal{T}_{m,k},\chi\rangle}\\
& \sim&
\frac{k}{2\pi}\cdot \int_X\int_{-\infty}^{+\infty}\int_{-\epsilon}^{+\epsilon}\chi(\lambda)\,\beta(\tau)\,e^{-ik\tau\lambda}
\,U(\tau)(x,y)\Pi_k(y,x)\,d\mu_X(y)\,d\lambda\,d\tau,\nonumber
\end{eqnarray}
where $\sim$ stands for \lq equal asymptotics for $k\rightarrow +\infty$\rq.
Remark that
\begin{eqnarray}
\label{eqn:fourier-of-k}
\Pi_k(y,x)&=&\frac{1}{2\pi}\int_{-\pi}^{\pi}e^{-ik\vartheta}\,\Pi\big(r_\vartheta(y),x\big)\,d\vartheta\nonumber\\
&=&\frac{1}{2\pi}\int_{-\pi}^{\pi}e^{-ik\vartheta}\,\Pi\big(y,r_{-\vartheta}(x)\big)\,d\vartheta,
\end{eqnarray}
where $r_\vartheta:X\rightarrow X$ is the action of $e^{i\vartheta}\in S^1$.

For any $s>0$, let $M_s(m)=:\left\{m'\in M:\mathrm{dist}_M\left(m',m\right)<s\right\}$,
where $\mathrm{dist}_M$ is the Riemannian distance of $M$.
Choose a sufficiently small $\delta>0$, and let $M_1=:M_{2\delta}(m)$, $M_2=:\overline{M_{\delta}(m)}^c$.
Let furthermore
$\{\varrho_1,\varrho_2\}$ be a partition of unity on $M$, subordinate to the open cover $\{M_1,M_2\}$, and
write $\varrho_j$ for $\varrho_j\circ \pi$.
Inserting the identity $\varrho_1(y)+\varrho_2(y)=1$ in (\ref{eqn:fourier-transform-comp}), we get
$\langle \mathcal{T}_{m,k},\chi\rangle \sim \langle \mathcal{T}_{m,k},\chi\rangle_1 +\langle \mathcal{T}_{m,k},\chi\rangle_2$, where
in $\langle\mathcal{T}_{m,k},\chi\rangle_j$ the integrand has been multiplied by $\varrho_j(y)$.

\begin{lem} As $k\rightarrow +\infty$, we have
$\langle\mathcal{T}_{m,k},\chi\rangle_2=O\left(k^{-\infty}\right)$.
\end{lem}

\textit{Proof.} Since the singular support of $\Pi$ is the diagonal in $X\times X$,
on the support of $ h\left(y,x'\right)=:\varrho_2(y)\,\Pi\left(y,x'\right)$ we have
$\mathrm{dist}_M\big(\pi(y),\pi(x)\big)\ge \delta$; therefore $h$
is $\mathcal{C}^\infty$ on $X_2\times X'$, where $X_2=\pi^{-1}(M_2)$ and $X'=\pi^{-1}\left(M'\right)$, with $M'$
a suitably small open neighborhood of $m$ (say, $M'=M_{\delta/2}(m)$).
We may thus regard $h\left(\cdot,x'\right)$ as a $\mathcal{C}^\infty$ family of compactly supported
$\mathcal{C}^\infty$-functions on
$X_2$, parameterized by $x'\in X'$, where both $X_2$ and $X'$ are $S^1$-invariant.

For each $\tau\in \mathbb{R}$, on the other hand, $U(\tau)$ is
a FIO associated to a canonical graph. Therefore, we may also regard $\kappa\left(\tau,x''\right)=:U(\tau)\left(x'',\cdot\right)$
as a $\mathcal{C}^\infty$ family of distributions on $X$, parametrized by $\left(\tau,x''\right)$.
Explicitly, $\left<\kappa\left(\tau,x''\right),f\right>=\big(U(\tau)(f)\big)\left(x''\right)$ for
$f\in \mathcal{C}^\infty(X)$.

Hence $B\left(\tau,x',x''\right)=:\langle \kappa\left(\tau,x''\right),h\left(\cdot,x'\right)\rangle$
is a $\mathcal{C}^\infty$-function of $\left(\tau,x',x''\right)\in \mathbb{R}\times X'\times X$, and therefore
its $l$-th Fourier component $B_l\left(\tau,x',x''\right)$
with respect to $x'$ is $O\left(l^{-\infty}\right)$ as $l\rightarrow \infty$,
uniformly on compact subsets of $\mathbb{R}\times X'\times X$.

Now in view of (\ref{eqn:fourier-transform-comp}), (\ref{eqn:fourier-of-k})
and the above, we have
\begin{eqnarray*}
\langle \mathcal{T}_{m,k},\chi\rangle_2&=&
\frac{k}{2\pi}\cdot \int_{-\infty}^{+\infty}\int_{-\epsilon}^{+\epsilon}\chi(\lambda)e^{-ik\tau\lambda}
\,\beta(\tau)\,B_{-k}(\tau,x,x)\,d\lambda\,d\tau\nonumber\\
&=&\frac{k}{\sqrt{2\pi}}\cdot \int_{-\epsilon}^{+\epsilon}\widehat{\chi}(k\tau)\,\beta(\tau)
\,B_{-k}(\tau,x,x)\,d\tau=O\left(k^{-\infty}\right).
\end{eqnarray*}

\hfill Q.E.D.

\bigskip

As in \cite{p-weyl}, we shall make use of the microlocal structure of $U(\tau)$ and $\Pi$:
up to smoothing terms that contribute negligibly to
the asymptotics, these may be represented as Fourier integral operators.

Namely, working in local coordinates near $x$ on the one hand we may write for $\tau\sim 0$:
\begin{eqnarray}
\label{eqn:fourier-for-U}
U(\tau)\left(x',x''\right)=\frac{1}{(2\pi)^{2\mathrm{d}+1}}\,\int_{\mathbb{R}^{2\mathrm{d}+1}}
e^{i[\varphi(\tau,x',\eta)-x''\cdot \eta]}a\left(\tau,x',x'',\eta\right)\,d\eta,
\end{eqnarray}
where the generating function $\varphi$ and the amplitude $a$ are as follows \cite{gs}.
First,
\begin{equation}
\label{eqn:fase-generatrice}
\varphi(\tau,x',\eta)=x'\cdot \eta+\tau\,q\left(x',\eta\right)+O\left(\tau^2\right)\,\|\eta\|;
\end{equation}
furthermore, $a(\tau,\cdot,\cdot,\cdot)\in S^0_{\mathrm{cl}}$ for every $\tau$, with
$a\left(0,x',x'',\eta\right)=1/\mathcal{V}\left(x''\right)$, where $\mathcal{V}$ is the local coordinate expression
of the volume form on $X$.

On the other hand, after \cite{bs}, $\Pi$ is a Fourier integral operator with complex phase:
\begin{equation}
\label{eqn:fourier-for-pi}
\Pi\left(x',x''\right)=\int_0^{+\infty}e^{i t\psi\left(x',x''\right)}s\left(t,x',x''\right)\,dt,
\end{equation}
where the Taylor expansion of the phase $\psi$ along the diagonal is determined by the K\"{a}hler metric, and
$s\left(t,x',x''\right)\sim \sum_{j\ge 0}t^{\mathrm{d}-j}s_j\left(x',x''\right)$.

Inserting (\ref{eqn:fourier-for-U}) and (\ref{eqn:fourier-for-pi}) in (\ref{eqn:fourier-transform-1}),
we get
\begin{eqnarray}
\label{eqn:insert-foi}
\lefteqn{\langle \mathcal{T}_{m,k},\chi\rangle\sim
\frac{k}{(2\pi)^{2\mathrm{d}+3}}
}\\
&&\cdot  \int_{-\infty}^{+\infty}\int_{-\epsilon}^{+\epsilon}\int_{\mathbb{R}^{2\mathrm{d}+1}}\int_{-\pi}^\pi\int_X
e^{i\Phi_k}\,H_k\big(\lambda,\tau,x,y,\eta,t,\vartheta\big)\,d\lambda\,d\tau\,d\eta\,d\vartheta\,d\mu_X(y), \nonumber
\end{eqnarray}
where
\begin{eqnarray*}
\lefteqn{\Phi_k=\Phi_k(t,\vartheta,\lambda,\tau,x,y,\eta)}\\
&=:&\varphi(\tau,x,\eta)-y\cdot \eta-k\lambda \tau+t\,\psi\big(r_\vartheta(y),x\big)-k\vartheta\end{eqnarray*}
$$
H_k\big(t,\vartheta,\lambda,\tau,x,y,\eta\big)=:\chi(\lambda)\,\beta(\tau)\,\varrho_1(y)\,
\,a\left(\tau,x,y,\eta\right)\,s\big(t,r_\vartheta(y),x\big).
$$

With the change of variables $\eta\rightarrow k\eta$ and $t\rightarrow kt$, (\ref{eqn:insert-foi}) may be rewritten
\begin{eqnarray}
\label{eqn:insert-foi-rescaled}
\lefteqn{\langle \mathcal{T}_{m,k},\chi\rangle\sim
\left(\frac{k}{2\pi}\right)^{2\mathrm{d}+3}}\\
&&\cdot  \int_{-\infty}^{+\infty}\int_{-\epsilon}^{+\epsilon}\int_{\mathbb{R}^{2\mathrm{d}+1}}\int_{0}^{+\infty}\int_{-\pi}^\pi\int_X
e^{ik\Phi}\,H_k\big(\lambda,\tau,x,y,k\eta,k t,\vartheta\big)\,d\lambda\,d\tau\,d\eta\,dt\,d\vartheta\,d\mu_X(y), \nonumber
\end{eqnarray}
where
$$
\Phi(t,\vartheta,\lambda,\tau,x,y,\eta)=\varphi(\tau,x,\eta)-y\cdot \eta-\lambda \tau+t\,\psi\big(r_\vartheta(y),x\big)-\vartheta.
$$

Implicitly introducing a partition of unity in $\vartheta$, we may assume that the integrand is compactly
supported in $\vartheta$. We now make the following remarks, that may proved by an adaptation of the arguments
in Lemmata 2.3 - 2.5 of \cite{p-weyl} (the focus in \cite{p-weyl} is on the asymptotics with respect to a continuous
parameter denoted $\lambda$, while here the asymptotics are with respect to the discrete parameter $k$).

\begin{itemize}
  \item Integration by parts in $d\vartheta$ shows that only a negligible contribution to the asymptotics is lost if
we multiply the integrand by a compactly supported function $\zeta=\zeta (t)$ with
$\zeta\in \mathcal{C}^\infty _0\big((1/C,C)\big)$ and $\zeta=1$ on $(2/C,C/2)$.
  \item We may thus replace $\int_{0}^{+\infty}dt$
by $\int_{0}^{+\infty}\zeta(t)\,dt$ in (\ref{eqn:insert-foi-rescaled}).
  \item Integration by parts in $dt$ then shows that for any $c>0$ and $\xi<1/2$
the contribution to the asymptotics coming from the loci $W_k\subseteq X$ where $\mathrm{dist}_M(y,x)\ge c\,k^{-\xi}$
is negligible.
\end{itemize}
Let us fix $\xi\in [1/3,1/2)$.

Next let us choose a system of Heisenberg local coordinates for $X$ centered at $x$
(we refer to \cite{sz} for a definition and discussion of Heisenberg local coordinates); let
$x+(\theta,\mathbf{v})$ denote the point with Heisenberg local coordinates $(\theta,\mathbf{v})$.
Thus we write $y=x+(\theta,\mathbf{v})$ and $d\mu_X(y)=\mathcal{V}(\theta,\mathbf{v})\,d\theta\,d\mathbf{v}$; in particular,
$\mathcal{V}(\theta,\mathbf{0})=1/(2\pi)$ for every $\theta$. By the previous discussion, only a negligible
contribution to the asymptotics is lost if we multiply the integrand by $\gamma\left(k^{\xi}\|\mathbf{v}\|\right)$, where
$\gamma\in \mathcal{C}^\infty_0(\mathbb{R})$ satisfies $\gamma(b)=1$ if $|b|\le 1$,
$\gamma(b)=0$ if $|b|\ge 2$.

Since this amounts to removing a smoothing term from $U(\tau)$, we may multiply the amplitude $a$ by a cut-off function
in $\eta$ which vanishes for $\|\eta\|<\delta$ for some small $\delta>0$, and equals $1$ for large $\|\eta\|$.
We may then also write $\eta=r\omega$, with $r>0$ and $\omega\in S^{2\mathrm{d}}$, so that $d\eta=r^{2\mathrm{d}}\,dr\,d\omega$,
and replace $\int_{\mathbb{R}^{2\mathrm{d}+1}}d\eta$ by $\int_0^{+\infty}\int_{S^{2\mathrm{d}}}r^{2\mathrm{d}}\,dr\,d\omega$.
We shall write $\omega=(\omega_0,\omega_1)\in S^{2\mathrm{d}}\subseteq\mathbb{R}\times \mathbb{R}^{2\mathrm{d}}$
with $\omega_0^2+\|\omega_1\|^2=1$.

Now $\big((0,\mathbf{0}),(1,\mathbf{0})\big)$
corresponds to $(x,\alpha_x)$ in Heisenberg local coordinates.
An argument similar to the proof of Lemma 2.2 of \cite{p-weyl} shows that as
$k\rightarrow +\infty$ only a rapidly decreasing
contribution is lost if integration in $d\omega$ is restricted
to a small open neighborhood $S_+\subseteq S^{2\mathrm{d}}$ of $(1,\mathbf{0})$;
we may thus multiply the integrand in (\ref{eqn:insert-foi-rescaled}) by an appropriate bump
function $g\in \mathcal{C}^\infty_0(S_+)$ which is identically $1$ near $(1,\mathbf{0})$
without affecting the asymptotics. We may assume $\omega_0\ge a>0$ for some fixed $a>0$
on $S_+$.

We shall furthermore adopt the rescaling
$\mathbf{v}\rightarrow \mathbf{v}/(r\sqrt{k})$, and write
$y=x+\big(\theta,\mathbf{v}/(r\sqrt{k})\big)$, so that
$d\mu_X(y)=r^{-2\mathrm{d}}k^{-\mathrm{d}}\mathcal{V}(\theta,\mathbf{v})\,d\theta\,d\mathbf{v}$ in the new coordinates. In rescaled coordinates,
integration in $d\mathbf{v}$ is over a ball in $\mathbb{C}^\mathrm{d}$ centered at the origin of radius $2\,r\,k^{1/2-\xi}$.

Given this, (\ref{eqn:insert-foi-rescaled}) may be rewritten

\begin{eqnarray}
\label{eqn:insert-foi-rescaled-in-v}
\lefteqn{\langle \mathcal{T}_{m,k},\chi\rangle\sim
\frac{k^{\mathrm{d}+3}}{(2\pi)^{2\mathrm{d}+3}}}\\
&&\cdot  \int_{-\infty}^{+\infty}\int_{-\epsilon}^{+\epsilon}\int_0^{+\infty}\int_{S^{2\mathrm{d}}}
\int_{1/C}^{C}\int_{-\pi}^\pi\int_{-\pi}^\pi\int_{\mathbb{C}^\mathrm{d}}
e^{ik\Psi_k}\,S_k\,d\lambda\,d\tau\,dr\,d\omega\,dt\,d\vartheta\,d\theta\,d\mathbf{v}, \nonumber
\end{eqnarray}
where
\begin{eqnarray*}
\lefteqn{\Psi_k\big(t,\vartheta,\lambda,\tau,x,\theta,\mathbf{v},r,\omega\big)=:
\Phi\left(\lambda,\tau,x,x+\left(\theta,\frac{\mathbf{v}}{r\sqrt{k}}\right),r\omega,\vartheta\right)}\\
&=&-\frac{1}{\sqrt{k}}\,\mathbf{v}\cdot \omega_1+\Big[-r\omega_0\theta+\tau r\,q(x,\omega)
+O\left(\tau^2\right)\,r\\
&&\left.-\lambda\,\tau+t\,\psi\left(x+\left(\theta+\vartheta,\frac{\mathbf{v}}{r\sqrt{k}}\right),x\right)-\vartheta\right],
\end{eqnarray*}
and
\begin{eqnarray*}
S_k\big(\lambda,\tau,x,\theta,\mathbf{v},r,\omega,\vartheta\big)=:H_k\big(\lambda,\tau,x,y,k\eta,k t,\vartheta\big)\,
\zeta(t)\,\gamma\left(k^{\xi-1/2}\|\mathbf{v}\|/r\right)\,g(\omega).
\end{eqnarray*}

By (65) of \cite{sz}, we have
\begin{eqnarray}
\label{eqn:exp-for-psi}
\lefteqn{t\,\psi\left(x+\left(\theta+\vartheta,\frac{\mathbf{v}}{r\sqrt{k}}\right),x\right)}\\
&=&i t\,\left[1-e^{i(\theta+\vartheta)}\right]+\left[i\,\frac{\|\mathbf{v}\|^2 }{2 r^2 k}
+R^\psi_3\left(\frac{\mathbf{v}}{r\sqrt{k}}\right)\right]\,t\,e^{i(\theta+\vartheta)},\nonumber
\end{eqnarray}
where $R^\psi_3$ vanishes to third order at the origin.

We can then write
\begin{eqnarray*}
ik\Psi_k&=&-i\,\sqrt{k}\,\mathbf{v}\cdot \omega_1+i\,k\,\Psi-\frac{\|\mathbf{v}\|^2 }{2 r^2 }\,t\,e^{i(\theta+\vartheta)}
+i\,k\,R^\psi_3\left(\frac{\mathbf{v}}{r\sqrt{k}}\right)\,t\,e^{i(\theta+\vartheta)},\nonumber
\end{eqnarray*}
where
\begin{eqnarray*}
\lefteqn{\Psi\big(t,\theta,\vartheta,\lambda,\tau,x,\mathbf{v},r,\omega\big)}\\
&=:&-r\omega_0\theta+\tau r\,q(x,\omega)
+O\left(\tau^2\right)\,r-\lambda\,\tau+i t\,\left[1-e^{i(\theta+\vartheta)}\right]-\vartheta.
\end{eqnarray*}

Since the exponential $D_k(\mathbf{v},\theta,\vartheta,t)=:\exp\left(i\,k\,R^\psi_3\big(\mathbf{v}/(r\sqrt{k})\big)\,t\,e^{i(\theta+\vartheta)}\right)$
is bounded in the range $\|\mathbf{v}\|<2\,k^{1/2-\xi}$, it may be incorporated into the amplitude (as in \S 3 and
\S 5 of  \cite{sz}). Furthermore, since $\partial_\theta\Psi=-r\omega_0+t\,e^{i(\theta+\vartheta)}$ and
$0<a\le\omega_0<1$, $1/C\le t\le C$, integration by parts in $\theta$ shows that, for some fixed $D\gg 0$,
the regions where $r<1/D$ or $r>D$ only give a negligible contribution to the
asymptotics for $k\rightarrow +\infty$.

In the same manner, integration by parts in $d\tau$ shows that only a negligible contribution is lost if
the integrand is multiplied by $\sigma(\lambda)$, where $\sigma\in \mathcal{C}^{\infty}_{0}\big((-C,C)\big)$
and $\sigma=1$ on $(-C/2,C/2)$.

Replacing $D$ and $C$ with $\max\{D,C\}$,
and multiplying the integrand by a suitable bump function $\kappa (r)$, we obtain
\begin{eqnarray}
\label{eqn:insert-foi-inner-outer}
\langle \mathcal{T}_{m,k},\chi\rangle\sim
\frac{k^{\mathrm{d}+3}}{(2\pi)^{2\mathrm{d}+3}}\int_{\mathbb{C}^\mathrm{d}}\int_{S^{2\mathrm{d}}}
e^{-i\sqrt{k}\,\mathbf{v}\cdot\omega_1}\,I_k(\mathbf{v},\omega)\,d\mathbf{v}d\omega,
\end{eqnarray}
where
\begin{equation}
\label{eqn:defn-di-I-k}
I_k(\mathbf{v},\omega)=:\int_{-C}^{C}\int_{-\epsilon}^{+\epsilon}\int_{1/C}^{C}
\int_{1/C}^{C}\int_{-\pi}^\pi\int_{-\pi}^\pi
e^{ik\Psi_{\mathbf{v},\omega}}\,\sigma(\lambda)\cdot T_k\,d\lambda\,d\tau\,dr\,dt\,d\vartheta\,d\theta\, \nonumber
\end{equation}
with
$T_k=: S_k\cdot D_k$ and
$$
\Psi_{\mathbf{v},\omega}\big(t,\theta,\vartheta,\lambda,\tau,r\big)=:
\Psi\big(t,\theta,\vartheta,\lambda,\tau,x,\mathbf{v},r,\omega\big).
$$
Thus we first view $\mathbf{v}$ and $\omega$ as parameters in the inner integral, and
evaluate it asymptotically using the stationary phase Lemma. A straightforward computation leads to
the following:

\begin{lem}
\label{lem:stationary-point-psi}
$\Psi_{\mathbf{v},\omega}$ has a unique stationary point
$$
\big(t_0,\theta_0,\vartheta_0,\lambda_0,\tau_0,r_0\big)=
\Big(1,0,0,q(x,\omega)/\omega_0,0,1/\omega_0\Big).
$$
The Hessian of $\Psi_{\mathbf{v},\omega}$ at the critical point has determinant
$$
\det\Big(H(\Psi_{\mathbf{v},\omega})\Big)=-\omega_0^2.
$$
\end{lem}

In particular, $\Psi_{\mathbf{v},\omega}$ vanishes at the critical point, and
$$
\frac{1}{\sqrt{\det\big(kH(\phi)/2\pi i\big)}}=\frac{1}{\omega_0}\,\left(\frac{k}{2\pi}\right)^{-3}.
$$
Applying the stationary phase Lemma, the asymptotic expansion of the amplitudes of $\Pi$ and
$U(\tau)$ and the Taylor expansion of $s_j$, $a_j$, $D_k$ and $\mathcal{V}$, we conclude that for every integer $N\gg 0$
\begin{eqnarray}
\label{eqn:stima-di-I-k}
\lefteqn{I_k(\mathbf{v},\omega)\sim \frac{1}{\omega_0}\,g(\omega)\,\gamma\left(k^{\xi-1/2}\omega_0\,\mathbf{v}\right)
e^{-\omega_0^2\,\|\mathbf{v}\|^2/2}}\\
&&\cdot \left(\frac{k}{2\pi}\right)^{-3}\,\left(\frac k\pi\right)^{\mathrm{d}}\,
\left[\chi\left(\frac{q(x,\omega)}{\omega_0}\right)+\sum_{j= 1}^N k^{-j/2} P_j(\mathbf{v},\omega,\chi)\left(\frac{q(x,\omega)}{\omega_0}\right)\right]
+R_N,\nonumber
\end{eqnarray}
where $\big|R_N\big|<C'_N\,k^{-aN}e^{-a\|\mathbf{v}\|^2}Q_N(\mathbf{v})$, for some $C'_N,a>0$ and some
polynomial $Q_N$, while each $P_j$ is a differential operator of degree $2j$ in $\chi$, as prescribed
by the stationary point Lemma, and is
a polynomial in $\mathbf{v}$;
we have made made use of the equalities $s_0(x,x)=\pi^{-\mathrm{d}}$, $a(0,x,y,\eta)=1/\mathcal{V}(y)$,
$\mathcal{V}(\theta,\mathbf{0})=1/2\pi$.

In the way to estimate the outer integral in (\ref{eqn:insert-foi-inner-outer}), let us first
remark that, for $N\gg 0$,
integration of the remainder $R_N$ over a ball of
radius $O\left(k^{1/2-\xi}\right)$ is $O\left(k^{-a'N}\right)$ for some $a'>0$. To estimate the integral of
the former summand, we may apply the stationary phase Lemma in $\mu=\sqrt{k}$.  The phase $\Upsilon=-\mathbf{v}\cdot \omega_1$ has a unique stationary
point for $\mathbf{v}=\omega_1=\mathbf{0}$, and partial integration in $d\omega$ shows that
only a negligible contribution
to the asymptotics is lost if integration in $d\mathbf{v}$ is restricted to an arbitrary open neighborhood of the origin; we
may thus replace $\gamma\left(k^{\xi-1/2}\omega_0\,\mathbf{v}\right)$ with some fixed cut-off $\rho=\rho(\mathbf{v})$
identically one near the origin. At the critical point $\Upsilon$ clearly vanishes, and
its Hessian $\Upsilon''$ satisfies
$$
\det\left(\frac{\mu \Upsilon''}{2\pi i}\right)=\left(\frac{\mu}{2\pi}\right)^{4\mathrm{d}}
=\frac{k^{2\mathrm{d}}}{(2\pi)^{4\mathrm{d}}}.
$$
Since $q\big(x,(1,\mathbf{0})\big)=q(x,\alpha_x)=\varsigma_T(m)$, we get
\begin{eqnarray}
\label{eqn:final-expansion}
\langle \mathcal{T}_{m,k},\chi\rangle\sim
\left(\frac k\pi\right)^{\mathrm{d}}\left[\chi\big(\varsigma_T(m)\big)+
\sum_{j\ge 1}k^{-j/2}Q_j(\chi)\big(\varsigma_T(m)\big)\right],
\end{eqnarray}
where the $Q_j$'s are differential operators acting on $\chi$.

Now we remark that the asymptotic expansions for the amplitudes $s$ and $a$ of $\Pi$ and $U$, respectively,
go down by integer steps; therefore, the appearance of fractional powers of $k$ in (\ref{eqn:stima-di-I-k})
is due to the asymptotic expansion of the amplitude in $\mathbf{v}/\sqrt{k}$. Hence the general summand in
(\ref{eqn:stima-di-I-k}) splits as a sum of multiples of $k^{d-3-s/2-l}Q_{s,l}(\mathbf{v},\omega,\chi)$, where
$s$ and $l$ are non-negative integers, and $Q_{s,l}(\mathbf{v},\omega,\chi)$
is a homogeneous polynomial of degree $s$ in $\mathbf{v}$, and a differential
operator of degree (at most) $2l$ in $\chi$. In turn, when we apply the stationary phase Lemma in $\mu$,
and keep track of powers of $k$ involved, this summand gives rise to a linear combination of terms of the form
$$
k^{\mathrm{d}+3}\, k^{-\mathrm{d}}\, k^{d-3-s/2-l}\,k^{-t/2}\Big(\frac{\partial}{\partial v_{i_1}}\circ
\frac{\partial}{\partial \omega_{1i_1}}\circ \cdots \circ \frac{\partial}{\partial v_{i_t}}\circ
\frac{\partial}{\partial \omega_{1i_t}}\Big)Q_{s,l}(\mathbf{0},(1,\mathbf{0}),\chi),
$$
with $t$ a non-negative integer. Since $Q_{s,l}$ is homogeneous of degree $s$ in $\mathbf{v}$ and we are evaluating at the
origin, we only get a non-zero contribution for $t=s$. Hence, we get a linear combination of summands of the
form
$$
k^{d-s-l}P_{s,l}(\chi)\big(\varsigma_T(m)\big),
$$
with $P_{s,l}$ a differential operator in $\chi$, of degree $2l\le 2(s+l)$.

\hfill Q.E.D.

\end{document}